\documentclass[11pt]{article}

\usepackage{latexsym}
\usepackage{amsmath}
\usepackage{amssymb}
\usepackage{dsfont}
\usepackage{amsmath}
\usepackage{fancyhdr}
\usepackage{mathrsfs}
\usepackage{times}
\usepackage{graphicx}
\usepackage{theorem}

\def \N {\mathbb N}
\def \C {\mathbb C}

\def \R {\mathbb R}

\def \F {\mathbb F}
\def \P {\mathbb P}
\def \O {\mathcal O}
\def \M {\mathcal M}
\def \W {\mathcal W}
\def \F {\mathcal F}
\def \NB {\mathcal N}
\def \B {\mathcal B}
\def \J {\mathcal J}
\def \H {\mathcal H}
\def \L {\mathcal L}
\def \U {\mathcal U}
\def \V {\mathcal V}
\begin{document}

\author{Song Sun}

\title{Note on geodesic rays and simple test configurations}
\date{}
\maketitle

 \begin{abstract}
In this short note, we give a new proof of a theorem of
Arezzo-Tian on the existence of smooth geodesic rays tamed by a
special degeneration with smooth central fiber.
 \end{abstract}
In \cite{D1}, S. K. Donaldson proposed a program to tackle
 the problem of the existence and uniqueness of extremal metrics on a K\"{a}hler manifold from the perspective of the
infinite dimensional space of K\"{a}hler potentials. He observed
that the existence of smooth geodesics connecting two arbitrary
K\"{a}hler potentials implies the uniqueness of K\"{a}hler metrics
in the given class with constant scalar curvature.  In \cite{C1},
X-X. Chen proved the existence of $C^{1,1}$ geodesics joining two
arbitrary points in $\H.\;$ Consequently, this established the
uniqueness of extremal K\"ahler metrics when the first Chern class
of the manifold is non-positive.   At present, there is extensive
research in this direction.
In particular, the uniqueness problem has been completely settled (c.f. \cite{M2},\cite{D4} and \cite{CT2}).\\

   We shall first give a very brief outline about a small part of this program
   which is directly relevant to the problem at hand.
   For more detailed accounts, readers are referred to \cite{D1}, \cite{C1},  \cite{CT2} and \cite{C2}.\\

Let $(M, \omega, J)$ be an $n$ dimensional K\"{a}hler manifold.
Define the infinite dimensional space of K\"{a}hler potentials as
 $$\H=\{\phi\in C^{\infty}(M)|\omega_{\phi}=\omega+\sqrt{-1}\partial\overline{\partial}\phi>0\}. $$
In \cite{M1}(c.f. \cite{D1}, \cite{S}), T. Mabuchi first
introduced a Weil-Petersson type metric  on $\H:$
\begin{equation*}
(\phi_1,\phi_2)_{\phi}=\int_{M}\phi_1\phi_2\frac{\omega_{\phi}^n}{n!},
\end{equation*}
 where $\phi_1, \phi_2\in T_{\phi}\H\simeq C^{\infty}(M)$.  It
is easy to see that the geodesic equation in $\H$  is\\
\begin{equation}
\ddot{\phi}=\frac{1}{2}|\nabla_{\phi}\dot{\phi}|_{\phi}^2
\end{equation}\\
A straightforward calculation shows (c.f. \cite{D1}, \cite{M1},
\cite{S}) that the space $\H$ is formally of non-positive
curvature. This fact was made rigorously in \cite{CC}, where  E.
Calabi and X-X. Chen proved that $\H$ is a non-positively curved
space
in the sense of Alexanderov.\\

According to S.Semmes \cite{S}, by adding a trivial $S^1$ factor,
the geodesic equation could be written as a degenerate complex
Monge-Amp\`{e}re equation in $M \times ([0,1] \times S^1).\;$
Suppose $X$ is a  Riemann surface with boundary. Denote  $\pi_1:
M\times X\rightarrow M$ and $\pi_2: M\times X\rightarrow X$ as the
two natural projection maps, and let $\Omega=\pi_1^{*}\omega.\;$
Then, given $\phi_0\in C^{\infty}(M\times\partial X)$ such that
$\Omega+\sqrt{-1}\partial\overline{\partial}\phi_0>0$ on each
slice
$M\times \{x\}$ for all $ x\in \partial X$, we consider the Dirichlet boundary value problem:\\
\begin{equation}
\left\{
                                                               \begin{array}{ll}
(\Omega+\sqrt{-1}\partial\overline{\partial}\Phi)^{n+1}=0,   & \hbox{ on $M \times X$ ;} \\
          \Phi=\phi_0 ,    & \hbox{ on $M\times\partial X$.}\\

                                                               \end{array}
                                                             \right.
\end{equation}\\
 A solution is of geometric interest if
$\Omega+\sqrt{-1}\partial\overline{\partial}\Phi>0$ when
restricted on each slice $M\times\{x\}$ for all $x\in X$. Since
the target manifold $\H$ is an infinitesimal symmetric space,  any
smooth solution of (2) can be re-interpreted  (c.f. \cite{D1}) as
a
 harmonic map from $X$ to $\H$ with prescribed boundary
map $\phi_0:
\partial X \rightarrow \H.$  Any geodesic segment connecting $\phi_1$ with  $\phi_2$
corresponds to an $S^1$ invariant solution of (2) with
$X=[0,1]\times S^1$ and $\phi_0(0, \tau)=\phi_1(\tau)$, $\phi_0(1,
\tau)=\phi_2(\tau).\;$  The notion of a geodesic ray is similar to
the finite dimensional case: a geodesic ray in $\H$ is a geodesic
segment which can be infinitely extended  in one direction. In
other words,  a geodesic ray corresponds
to an $S^1$ invariant solution of the following:\\
\begin{equation}
(\Omega+\sqrt{-1}\partial\overline{\partial}\Phi)^{n+1}=0,\ \
\mbox{on} \ \ M \times ([0, \infty)\times S^1)\simeq M\times
(D\setminus\{0\}) ,
\end{equation}\\
   where $D$ is the closed
unit disk. \\

In \cite{D1}, Donaldson also  conjectured that  the existence of
smooth geodesic rays where the K energy is strictly decreasing at
the infinity is
 equivalent to the non-existence of constant scalar
curvature metrics in $[\omega]$. Donaldson's conjecture certainly
motivated the study of the existence of geodesic rays and related
problems.  However, the existence of geodesic rays is quite
different from that of geodesic segments since the domain involved
is naturally non-compact.    More importantly,  Donaldson
(\cite{D1}) pointed out that the initial value problem for the
geodesic ray equation is not always solvable in the smooth
category. So we need to impose an alternative condition in order
to properly solve equation (3).    Following Donaldson's program
\cite{D1}, this issue was discussed in \cite{C2}. According to
\cite{C2}, the initial K\"ahler potential together with the
asymptotic direction (given by either an existing geodesic ray or
an algebraic ray associated to a test configuration) forms a
well-posed Dirichlet boundary value for equation (3).   A set of
new problems were discussed there for developing the existence
theory of geodesic rays. In particular, it is proved the existence
of relative $C^{1,1}$ geodesic rays parallel to a given smooth
geodesic ray under natural geometrical constraints. Unfortunately,
there are few examples of the existence of geodesic rays in the
literature, which creates serious problem for pushing the general
existence theory further. In \cite{AT}, using Cauchy-Kowalewski
theorem, C. Arezzo and G. Tian proved the existence of a smooth
geodesic ray asymptotically parallel to a special degeneration
with smooth central fiber, or equivalently, to a simple test
configuration(c.f. \cite{D3}, \cite{CT1}).  One would like to see
a more direct PDE proof of this important theorem. The main
purpose of this note is to reprove the same theorem using
perturbation argument. Note that if we only assume the total space
of the test configuration to be smooth, then by \cite{CT1}, it
will give rise to a relative $C^{1,1}$ geodesic ray parallel to
the algebraic ray induced by the degeneration. Furthermore, it was
shown in \cite{CT1} that even for tori varieties, the best
regularity one could get is $C^{1,1}$(see also \cite{SZ}). There
are also other ways of obtaining geodesic rays
from a test configuration, see \cite{PS1}, \cite{PS2}.\\

Now we introduce the definition of a {\it K\"ahler fibration} and a {\it simple test configuration}.\\

\textbf{Definition 1}. A \emph{K\"{a}hler fibration} (over the
closed unit disk) is a map $\pi:(\M,J,\Omega)\rightarrow D$, where
$J$ is an integrable complex structure on $\M$, $\pi$ is a
holomorphic submersion, $\Omega$ is a closed two form on $\M$
which is compatible with $J$ and it is a K\"{a}hler form on each
fiber
$M_z(z\in D)$(which is assumed to be compact without boundary).\\

\textbf{Definition 2} (c.f. \cite{CT1}, \cite{D3}).  A (truncated)
\emph{simple test configuration} for a polarized K\"{a}hler
manifold $L\rightarrow M$ is a K\"{a}hler fibration $\pi:(\M, J,
\Omega)\rightarrow D$ together with a very ample line bundle $\L$
and a $\C^{*}$ equivariant embedding $\{\L\rightarrow\M\rightarrow
D\} \hookrightarrow \{\O(1)\rightarrow\P^{N}\times \C\rightarrow
\C\}$,
 such that $\{L\rightarrow M\}$ is isomorphic
$\{\L|_{M_1}\rightarrow M_1\}$, where we denote $M_t=\pi^{-1}(t)$.
Also, the $\C^{*}$ action on $\C$ is given by the standard
multiplication, and the map $\P^N\times\C\rightarrow \C$ is simply
the projection to the second factor, In addition, $\Omega$ should
coincide with the restriction of the Fubini-Study metric on
$\P^N$, while the induced $S^1$ actions on all these spaces are
assumed to be unitary. Clearly all the fibers $\pi^{-1}(t)$ for
$t\neq0$ are biholomorphic to each other. A simple test
configuration is called \emph{product} if $\M$ is biholomorphic to
$M\times \C$, and the $\C^{*}$ action on $\M$ is also a product
action coming from $\C^{*}$ action on $M$ and the standard
multiplication on $\C$. It is called \emph{trivial} if
the $\C^{*}$ action on $\M$ is also trivial.\\

\textbf{Remark 3}. The above definition of a simple test
configuration is essentially the same as the special degeneration
with smooth central fiber
studied  by G. Tian first in \cite{T}.\\

\textbf{Theorem 4} ( Arezzo-Tian  \cite{AT}). Given a non-trivial
simple test configuration for $L\rightarrow M$, there exists a non
trivial
geodesic ray tamed by this test configuration.\\

According to \cite{C2}, a geodesic ray is said to be {\it tamed by
a test configuration} if  it is asymptotically parallel to the
algebraic ray defined by pulling back the K\"{a}hler potentials
through the
$\C^{*}$ action on $\M$.\\

We want to take a different route to prove this theorem. Following
\cite{D2} and \cite{CT1},  smooth regular solutions to (3) are
related to foliations of punctured holomorphic discs with some
control on the total area.  There is a Fredholm theory associated
to the moduli space of holomorphic discs with totally real
boundary condition. Deformation of this moduli space is the
central topic of this note.
 \\

Arezzo-Tian's theorem  is a consequence of the following proposition.\\

 \textbf{Proposition 5}. let $\pi:(\M,J,\Omega)\rightarrow D$ be a K\"{a}hler
fibration, there exists a smooth function $\Phi$ defined in a
neighborhood of the central fiber $M_0$ that solves the complex
Monge-Amp\`{e}re equation
$(\Omega+\sqrt{-1}\partial\bar{\partial}\Phi)^{n+1}=0$ with $
\Omega+\sqrt{-1}\partial\bar{\partial}\Phi$ being positive on each fiber.\\

In \cite{AT}, it was shown that the value of $\Phi$ on the central
fiber could be prescribed as long as the corresponding K\"ahler
metric is real analytic. From the proof of proposition 5 we shall
see that also more is true. We only state it in the following
$S^1$ invariant case, since this
gives rise to geodesic rays. Namely,\\

\textbf{Theorem 6}. Let $\pi:(\M,J,\Omega)\rightarrow D$ be a
nontrivial simple test configuration. Denote $M_z=\pi^{-1}(z)$ for
$z\in D$, and $\omega_z=\Omega|_{M_z}$. For $k>0$ sufficiently
large, denote $\tilde{\H}_0=\{\phi\in C^{k+1}(M_0;\R)|\phi \ \
\text{is} \ \ S^1 \ \ \text{invariant and}\ \
\omega_0+\sqrt{-1}\partial\bar{\partial}\phi>0\}$. Given any
$\phi_0\in\tilde{\H}_0$, there is an open set $\U$ in
$\tilde{\H}_0$ and a number $\delta(k)>0$, such that for every
$\phi\in \U$ there exists a $S^1$ invariant $\Phi\in
C^k(\pi^{-1}(|z|\leq \delta))$ depending on $\phi$ that solves the
complex Monge-Amp\`{e}re equation
$(\Omega+\sqrt{-1}\partial\bar{\partial}\Phi)^{n+1}=0$ with $
\Omega+\sqrt{-1}\partial\bar{\partial}\Phi$ being positive on each
fiber.\\


The proof of Proposition 5 is based on a perturbation theory first
introduced in \cite{D2}  by Donaldson in the case of a trivial
test configuration. In this note, we follow its generalization in
\cite{CT1}. By the definition of a K\"{a}hler fibration, $\M$ is
always diffeomorphic to the product $M_0\times D.$  So we can for
simplicity assume $\M=M\times D$ for a $2n$ dimensional smooth
manifold $M$ and the map $\pi$ involved in the definition  is the
projection map to the second factor . Fix once and for all a cover
of $M\times D$ by small balls, say $\{U_i\}_{i\in I}$. Following
Donaldson's construction, we can associate a manifold $\W$ to any
K\"{a}hler fibration, as follows: On each $U_i$, we choose local
holomorphic coordinates to be $(z_1,\dots, z_n, z)$, where $z$ is
simply given by $\pi$. Then $\Omega$ could be written as
$\sqrt{-1}\partial\bar{\partial}\rho_i$ for some locally defined
function $\rho_i$. $\W$ is obtained by twisting the vertical
holomorphic cotangent bundle $E=T^{*}(M\times D)/\pi^{*}T^{*}D.\;$
More precisely, we glue $\xi$ in $E|_{U_i}$ with
$\xi+\partial(\rho_i-\rho_j)$ in $E|_{U_j}$ over the corresponding
fiber. It is easy to see that $\W$ is also a fibration over $D$
and the canonical complex-symplectic structure on the holomorphic
cotangent bundle induces a fiberwise complex-symplectic form on
$\W$. Furthermore,  $\Omega$ defines an exact LS-graph\footnote{In
a complex symplectic manifold $(M, \Theta)$, a submanifold $L$ is
called an \emph{LS-submanifold} if $L$ is Lagrangian with respect
to $Re\Theta$, while the restriction of $Im\Theta$
on $L$ is a symplectic form. For more details, see \cite{D2}, \cite{CT1}.} on each vertical fiber. \\

    Of course, our construction of $\W$ is not canonical. However, if we fix an open cover and an initial
K\"{a}hler fibration, then $\rho_i$ could be chosen to depend
smoothly on the data $\Omega$ and $J$ for a small
variation(Indeed, by the well known theorem of
Newlander-Nirenberg, holomorphic coordinates could be made to vary
smoothly. Then, one can follow the proof of Dolbeault's lemma to
show this). Moreover, by definition, $\W$ is always diffeomorphic
to $E$, or further, to the real vertical cotangent bundle, still
denoted by $E$, which is independent of $\Omega$ and $J$.
Therefore, if we pull back everything to the latter, a
perturbation of $\Omega$ and $J$ really gives us a perturbation of
the
complex-symplectic structure on $E$.\\

Now let $\phi_0: \partial D\rightarrow \R$ be a smooth function
such that $\Omega+\sqrt{-1}\partial\bar{\partial}\phi_0$ is
positive on fibers over $\partial D$. Then it defines exact
LS-graphs $\Lambda_{z,\phi_0}$ over any $z\in\partial D$.
Following \cite{D2}, \cite{CT1},
we have a one-to-one correspondence: \\

(A) A $C^{\infty}$ solution $\Phi$ to the homogeneous
Monge-Amp\`{e}re equation:
$(\Omega+\sqrt{-1}\partial\bar{\partial}\Phi)^{n+1}=0$ satisfying
the boundary condition $\Phi|_{\partial D}=\phi_0$ and such that
$\Omega+\sqrt{-1}\partial\bar{\partial}\Phi$ still
defines a K\"{a}hler fibration(together with $J$).\\

(B) A smooth map $G:M\times D\rightarrow E$ which covers the
identity map on $D$, holomorphic in the second variable, and
satisfies the boundary condition: for all $z\in \partial D$,
$G(\cdot, z)\in\Lambda_{z, \phi_0}$ (Alternatively, we could view
this as a family of holomorphic sections of the fibration
$E\rightarrow D$ whose boundary lies in some totally real
submanifold given by $\bigcup_{z\in \partial D}\Lambda_{z,
\phi_0}$). In addition, we require that $p_1\circ G(\cdot, 0)$ is
the identity map, and $p_1\circ G(\cdot, z)$ is a diffeomorphism
for any $z\in D$, where $p_1:E\rightarrow M$ is the
projection map.\\

\textbf{Lemma 7}. Perturbation of $\Omega$, $J$ and $\phi_0$
preserves a smooth solution to the above equation, i.e. the
compact family of normalized holomorphic discs in $(B)$ is stable
under
perturbation.\\

To prove this Lemma, we need to set up a Fredholm theory for
holomorphic discs with totally real boundary conditions. Denote by
$(D, E)_s$ the space of maps from $D$ to $E$ which lies in the
Sobolev space $H^{s+1}$ for some large $s$. Let $\F_s$ be the
subspace of $(D, E)_s$ which are sections of the fibration, i.e
normalized maps. Fix $J_0$ on $E$, and a totally real submanifold
$R_0$ of $E$ with respect to $J_0$(For example, in our case, the
exact Lagrangian graphs defined by the known smooth solution
restricted on $\partial D$). Denote by $NR_0$ the normal bundle of
$R_0$ in $E$ with respect to any fixed metric, then a $C^{\infty}
$ neighborhood of totally real submanifolds around $R_0$ can be
identified with an open set of the space $\Gamma(NR_0)$ of all
smooth sections of $NR_0$. As in \cite{O}, choose
$\bar{\epsilon}=\{\epsilon_k\}_{k\in\N}$($\epsilon_k\rightarrow
0$), and define a Floer norm on $\Gamma(NR_0)$ by
$$\parallel X\parallel_{\bar{\epsilon}}=\sum_{k\in\N}\epsilon_k \max_{x\in E}|D^k X(x)|.$$
For $r>0$, we also define
$$\Gamma^{\bar{\epsilon}}_r(NR_0)=\{X\in\Gamma(NR_0)|\parallel X\parallel_{\bar{\epsilon}}<r\}.$$
This is a Banach space and by choosing $r$ sufficiently small, we
can assume that $\Gamma^{\bar{\epsilon}}_r(NR_0)$ maps injectively
to the space of all totally real submanifolds of $E$, under the
previous identification. Let $\NB(R_0)$ denote its image. For each
$R\in\NB(R_0)$, there is an associated diffeomorphism $\phi_R:
R\rightarrow R_0$ which extends to a diffeomorphism of $E$.
Moreover, We can choose $\phi_R$ to depend smoothly on $R$. Now
let $\B=\cup_{u\in \F_s}H^s(u^{*}TE)$ be an infinite dimensional
vector bundle over $\F_s$, and $\J$ be the space of  smooth almost
complex structures on $E$ which are $C^{\bar{\epsilon}}$ close to
$J_0$. Then $\B\times(\partial D, E)_{s-\frac{1}{2}}$ is a bundle
over $\F_s\times \J\times\NB(R_0)$, with a section
$s(u,J,R)=(\bar{\partial_{J}}u, \phi_R^{-1}\circ u|_{\partial
D})$. Fix $J_0$, and let $s_0$ be the restriction of $s$ to the
slice $\F_s\times \{J_0\}\times \NB(R_0)$. It is shown in \cite{O}
that $s_0$ is transversal to the submanifold $\{0\}\times(\partial
D, R_0)_{s-\frac{1}{2}}$ at a point $(u_0, R_0)$ if $u_0$ is not
multiply covered, i.e. there exists a $z\in\partial D$, such that
$u_0^{-1}(u_0(z))\cap\partial D={z}$ and $Du_0(z)\neq0$. So in our
particular case $s$ is transversal to $\{0\}\times(\partial D,
R_0)_{s-\frac{1}{2}}$ at $(u_0, J_0, R_0)$ for every disc coming
from a solution of our previous equation (A). Therefore,
$s^{-1}(\{0\}\times(\partial D, R_0)_{s-\frac{1}{2}})$ is smooth
Banach manifold near $(u_0, J_0, R_0)$. Standard elliptic
regularity implies that it is actually contained in
$\F_{\infty}\times \J\times\NB(R_0)$ for $s$
sufficiently big.  \\

Now consider the projection map $s^{-1}(\{0\}\times(\partial D,
R_0)_{s-\frac{1}{2}})\rightarrow \J\times \NB(R_0)$, which is
Fredholm of index $2n$(c.f \cite{D2}, \cite{CT1}). Given a smooth
solution on $(M\times D, J, \Omega)$ as in $(A)$, we have a 2n
dimensional compact family of normalized holomorphic discs into
$(E, J_0)$, where $J_0$ is defined by $J$ and $\Omega$. Moreover,
the holomorphic discs appearing in the family are all
super-regular \footnote{For a family of holomorphic discs $G:
M\times D\rightarrow \W$ parameterized by $M$, we say that a disc
$G_x$($x\in M$) is super-regular if the derivative $dp_1\circ d_x
G(\cdot, z): T_xM\rightarrow T_{p_1\circ G(x, z)}M$ is surjective
for all $z\in D$. It is proved in \cite{D2}, \cite{CT1} that a
super-regular disc is automatically regular.}, and in particular
regular. Now if we perturb $J$, $\Omega$ and $\phi_0$, we are
actually perturbing $J_0$ and $R_0$. Standard Fredholm theory
ensures the existence of a nearby family of
normalized regular holomorphic discs, which proves Lemma 7. $\square$\\

\emph{Proof of proposition 5}. For $r\in(0,1)$, let $\M(r)$ be the
re-scaled K\"{a}hler fibration defined by $(\M, J,
\Omega)|_{|z|\leq r}$ with $\pi_r(w)=\pi(w)/r$. When $r$ is small
enough, $\M(r)$ is close to the trivial fibration given by the
product $(M_0, J|_{M_0}, \Omega|_{M_0})\times D$. The latter has
an obvious solution to $(A)$(just take $\Phi=0$). Therefore by
Lemma 7, for $r$ small, we obtain a solution to the equation on
$\M(r)$, which is the same as a solution near the central
fiber on $\M$. $\square$\\

\emph{Proof of Theorem 4}. The limit of the re-scaled test
configurations $\M(r)$ is $(M_0, J|_{M_0}, \Omega|_{M_0})\times D$
with a (possibly non-trivial) $C^*$ action on $M_0$. Any $S^1$
invariant potential $\phi_0$ on $M_0$ yields a trivial solution.
We can perturb $\phi_0$ to an $S^1$ invariant function $\phi$ on
$\partial\M(r)$ for small $r$, then the solution $\Phi$ ensured by
proposition 5 will also be $S^1$ invariant by the uniqueness of
solutions of equation (2), which follows from a standard maximum
principle argument(see lemma 6 in \cite{D1}). Then we obtain a
geodesic ray on the fiber $M_1$ by pulling back the restriction of
$\Omega+\sqrt{-1}\partial\bar{\partial}\Phi$ on each fiber to a
fixed fiber by the $\C^{*}$ action, and we also get a foliation by
punctured holomorphic discs on $M_1\times (D\setminus \{0\})$.
Furthermore, if the test configuration is non trivial, the
corresponding foliation would not be trivial since the $\C^{*}$
action on $\M$ is not along the leaf direction given by the
orthogonal complement of the tangent space of the fibers with
respect to $\Omega+\sqrt{-1}\partial\bar{\partial}\Phi$. Thus, in
this case, we do get a nontrivial geodesic ray. Since $\Phi$ is
smooth on $\M$, the geodesic ray is parallel to the algebraic ray
defined simply by pulling back $\Omega$ through the $\C^{*}$
action. Thus the ray is tamed by the
test configuration. $\square$\\

\emph{Proof of theorem 6}. First it is easy to see that the above
arguments still go through if we change the general framework
replacing $C^{\bar{\epsilon}}$ by $C^k$ for $k$ large. First of
all, following the proof of \cite{D2},  the one-to-one
correspondence before Lemma 7 becomes that a $C^{k+1}$ solution
$\Phi$ corresponds to a map $G$ which is only $C^k$ along $M$(i.e.
the corresponding exact LS-graph is only of class $C^k$). Now we
consider $C^k$ neighborhood $\NB^k(R_0)$ of a $C^k$ totally real
submanifold $R_0$. This could be identified with an open set of
the space $\Gamma^k(NR_0)$ of all $C^k$ sections of $NR_0$(in our
case, $R_0$ will be the exact LS-graphs defined over $\partial D$
corresponding to the boundary value of a $C^{k+1}$ solution to the
homogeneous Monge-Amp\`ere equation). More precisely, define
$$\Gamma^k_r(NR_0)=\{X\in \Gamma^k(NR_0)|\parallel X\parallel_{C^k}<r\}.$$
We also choose $r$ sufficiently small, and define $\NB^k(R_0)$ to
be the image of $\Gamma^k_r(NR_0)$ under the above identification.
Define $\F_s$, $\B$ and $\J$ the same as before, where $s$ should
be no bigger than $k$. The previous section $s(u, J,
R)=(\bar{\partial}_Ju, \phi_R^{-1}\circ u|_{\partial D})$ is $C^k$
in $u$, and $C^{\infty}$ in the remaining variables. The
submanifold of $ (\partial D, E)_{s-\frac{1}{2}}$ $$ (\partial D,
R_0)_{s-\frac{1}{2}}:=(\partial D, E)_{s-\frac{1}{2}}\cap
C^0(\partial D, R_0)$$ is a $C^k$ Banach manifold. So the implicit
function theorem implies that $s^{-1}(\{0\}\times (\partial D,
R_0)_{s-\frac{1}{2}})$ is $C^k$ Banach manifold near $(u_0, J_0,
R_0)$. Then the projection map $\pi_2: s^{-1}(\{0\}\times
(\partial D, R_0)_{s-\frac{1}{2}})\rightarrow \J\times \NB^k(R_0)
$ is a $C^k$ Fredholm map of index $2n$. $\pi_2^{-1}(J_0, R_0)$ is
a $C^k$ manifold near $(u_0, J_0, R_0)$ since it is parametrized
by $M$ under the map $G$. Since these consists of regular
holomorphic discs, we can find a neighborhood $\V$ of $(J_0,
R_0)$, such that there exists a $C^k$ diffeomorphism $F: \V\times
M\rightarrow \pi_2^{-1}(\V)$. Furthermore, since $s$ is smooth in
$\J$ and $\NB^k(R_0)$ direction, $F$ could be chosen to depend
smoothly on $\V$-variable. By using once again the correspondence
before lemma 7, we arrive at that $C^{k+1}$ perturbation of the
boundary condition gives rise to a $C^{k+1}$ perturbation of the
solution to (A), and the dependence is smooth.\\

With this at hand, we can define a smooth map $S$ from the open
set of elements $(r, \phi)$ in  $[0, 1]\times C^{k+1}(M;\R)$ such
that $\phi$ is the boundary value of a $C^{k+1}$ solution of the
homogeneous Monge-Amp\`ere equation for $\M(r)$ to $\tilde{\H}_0$.
Given $r$ and $\phi\in C^{k+1}(M;\R)$ in this open set, we take
the solution $\Phi$, and then restricts to the central fiber to
obtain an element in $\tilde{\H}_0$. For the product test
configuration $(M_0, J|M_0, \Omega|M_0)\times D$ with the $\C^*$
action on $M_0$, any $\phi\in\tilde{\H}_0$ forms a trivial
solution, and $S(\phi)=\phi$. The derivative with respect to the
second variable $(D_{2}S)_{(0, \phi)}: C^{k+1}(M;\R)\rightarrow
\tilde{\H}_0$ is then clearly surjective, with an obvious right
inverse which is the inclusion map. By the implicit function
theorem together with the previous paragraph, we conclude that for
any $\phi\in \tilde{\H}_0$, there exists a neighborhood $\U$ of
$\phi$, an $r>0$, and a smooth map $T: \U\rightarrow
C^{k+1}(\M(r); \R)$, such that the image of $T$
consists of solutions to (A) on $\M(r)$. $\square$\\

%

\textbf{Remark 8}. An interesting question is: Given a sequence of
K\"ahler potentials in $\H$ which is bounded in the sense of
Cheeger-Gromov, but not bounded in the holomorphic category. Does
there exist a point in the ``sphere at infinity" which reflects
this
non-compactness or degeneracy?\\

{\bf Acknowledgment}:  The author wishes to thank his advisor
Professor X-X. Chen for suggesting this problem and for kindly
turning down his request of co-authorship. He is also grateful to
the referee for helpful comments on the original version of this
article.

Department of Mathematics, University of Wisconsin-Madison, 480
Lincoln Dr, Madison, WI 53706.\\ Email:
ssun@math.wisc.edu\\\\

\end{document}